\documentclass[letterpaper, 10pt, conference]{ieeeconf}

\IEEEoverridecommandlockouts
\usepackage{amsmath, amsfonts, amssymb}
\usepackage{scalerel}

\renewcommand{\epsilon}{\varepsilon}
\renewcommand{\phi}{\varphi}

\newcommand{\kdmin}{\delta_{\min}}
\newcommand{\kdmax}{\delta_{\max}}
\newcommand{\numin}{\nu_{\min}}
\newcommand{\numax}{\nu_{\max}}

\newcommand{\mmin}{m_{\min}}
\newcommand{\mmax}{m_{\max}}
\newcommand{\sub}{^\mathrm{sub}}
\newcommand{\Loneloc}{L^1_\mathrm{loc}}

\newcommand{\dlt}{\scaleobj{0.7}{\triangle}}

\newcommand{\dif}{\mathrm{d}}

\newcommand{\calA}{\mathcal{A}}

\newcommand{\bbR}{\mathbb{R}}
\newcommand{\bbN}{\mathbb{N}}
\newcommand{\bfi}{\mathbf{i}}
\newcommand{\bfj}{\mathbf{j}}

\newcommand{\pb}[1]{\big( #1 \big)}

\newcommand{\px}[1]{\left( #1 \right)}

\newcommand{\abs}[1]{\vert #1 \vert}
\newcommand{\absb}[1]{\big\vert #1 \big\vert}

\newcommand{\setb}[1]{\big\lbrace #1 \big\rbrace}
\newcommand{\setB}[1]{\Big\lbrace #1 \Big\rbrace}

\newcommand{\norm}[1]{\Vert #1 \Vert}

\newtheorem{remark}{Remark}
\newtheorem{definition}{Definition}
\newtheorem{assumption}{Assumption}
\newtheorem{proposition}{Proposition}
\newtheorem{lemma}{Lemma}
\newtheorem{theorem}{Theorem}

\title{%\LARGE \bf
Stability of Two-Dimensional SISO LTI System with\\Bounded Feedback Gain That Has Bounded Derivative*
}

\author{Anton Ponomarev$^{1}$ and Lutz Gröll$^{1}$% <-this % stops a space
\thanks{*The authors gratefully acknowledge funding by the German Federal Ministry of Education and Research (BMBF) within the Kopernikus Project ENSURE ``New ENergy grid StructURes for the German Energiewende'' (03SFK1B0-3).}% <-this % stops a space
\thanks{$^{1}$The authors are with the Institute for Automation and Applied Informatics, Karlsruhe Institute of Technology, Eggenstein-Leopoldshafen, 76344 Baden-Württemberg, Germany {\tt\small \{anton.ponomarev, lutz.groell\}@kit.edu}}%
}

\begin{document}

\maketitle

\begin{abstract}
We consider a two-dimensional SISO LTI system closed by uncertain linear feedback. The feedback gain is time-varying, bounded, and has a bounded derivative (both bounds are known). We investigate the asymptotic stability of this system under all admissible behaviors of the gain. Note that the situation is similar to the classical absolute stability problem of Lurie--Aizerman with two differences: linearity and derivative constraint. Our method of analysis is therefore inspired by the variational ideas of Pyatnitskii, Barabanov, Margaliot, and others developed for the absolute stability problem. We derive the Hamilton--Jacobi--Bellman equation for a function describing the ``most unstable'' of the possible portraits of the closed-loop system. A numerical method is proposed for solving the equation. Based on the solution, sufficient conditions are formulated for the asymptotic stability and instability. The method is applied to an equation arising from the analysis of a power electronics synchronization circuit.
\end{abstract}

\section{Introduction}

The subject of the present paper is the single-input single-output (SISO) linear time-invariant (LTI) system
\begin{align*}
    \dot x &= Ax + bu, \\
    y &= c^T x
\end{align*}
closed by the linear time-varying feedback $u = \kappa(t) y$ whose gain $\kappa$ is unknown. What we do know are bounds on $\kappa$ and $\dot\kappa$. We investigate the asymptotic stability of the closed loop
\begin{equation*}
    \dot x = \pb{A + \kappa(t) B} x
\end{equation*}
where $B = bc^T$. It would be perfectly natural to treat the closed-loop dynamics as a \emph{linear parameter-varying} (LPV) system~\cite{shammaOverviewLPVSystems2012}. Popular results in the LPV theory are obtained using, e.g., quadratic Lyapunov functions~\cite{ducStabilityLinearTimevarying2006}, linear matrix inequalities~\cite{mozelliComputationalIssuesStability2019}, and integral quadratic constraints~\cite{jonssonSystemsUncertainParameters1996, pfiferLessConservativeRobustness2016}. Our approach, however, is rooted in the literature on \emph{absolute stability}. Let us explain relevant results from there.

The problem of absolute stability is one of the classical problems of control theory~\cite{liberzonEssaysAbsoluteStability2006, fradkovEarlyIdeasAbsolute2020}. For a SISO LTI system similar to the one above, absolute stability refers to the global asymptotic stability under an arbitrary \emph{nonlinear} feedback law $u = \phi(y,t)$~-- not a linear feedback like in our setting. Systems of this kind are called \emph{Lurie systems}. Absolute stability is commonly studied subject to Aizerman's \emph{sector condition} which, after appropriate scaling, reads
\begin{equation*}
    0 \leq \frac{\phi(y,t)}{y} \leq 1, \quad y \neq 0.
\end{equation*}
For our linear feedback, the sector condition translates to $\kappa(t) \in [0,1]$~-- an assumption we do make.

Some of the earliest absolute stability criteria are the \emph{circle criterion} and \emph{Popov criterion}~\cite{narendraFrequencyDomainCriteria1973}, the latter only applicable to the case of a time-invariant $\phi$. These conditions are formulated in the frequency domain. In a time-domain interpretation, they guarantee the existence of a time-invariant Lyapunov function structured as a quadratic form or an integral of the nonlinearity. Due to the presupposed structure of the Lyapunov function, the circle and Popov criteria are \emph{sufficient but not necessary}.

\emph{Necessary and sufficient} conditions in the absolute stability problem with a time-varying nonlinearity satisfying the sector condition can be attained by the \emph{variational method} proposed in~\cite{pyatnitskiiAbsoluteStabilityNonstationary1970}. An overview of the method is published in \cite{margaliotStabilityAnalysisSwitched2006}. The basic idea is to look at the ``worst-case'' or ``most destabilizing'' feedback~-- one with a switching gain that takes the system as far from the origin as possible. Then absolute stability is equivalent to the stability under the ``most destabilizing'' switching. For systems of \emph{dimension at most 3} this method yields practical necessary and sufficient conditions that can be checked numerically.

For \emph{2-dimensional} systems, the necessary and sufficient conditions have been known since the 1960s \cite{levinStabilitySolutionsSecondorder1961}. They can be derived from the geometric idea of a time-invariant \emph{comparison system}~\cite{leonovNecessarySufficientConditions2005} which is an event-triggered switching system designed so that its trajectories are \emph{transversal} to the trajectories of the closed-loop system, i.e., always crossed by the latter in a specified direction. The transversal trajectories are used to analyze the phase portrait of the closed loop in the spirit of the classical Poincaré--Bendixson analysis. The comparison system coincides with the ``most destabilizing'' switching system from the variational approach.

The variational/comparison ideas can lead to \emph{Lyapunov functions} that quantify absolute stability. A Lyapunov function that decreases at the exponential rate of the absolute (``worst-case'') Lyapunov exponent was defined and investigated in~\cite{barabanovAbsoluteCharacteristicExponent1989} using the variational concepts. This special Lyapunov function is now called the \emph{Barabanov norm}. A numerical method of its construction in the 2-dimensional case is developed in \cite{musaevaConstructionInvariantLyapunov2023}. The unit sphere of the Barabanov norm is piecewise analytical. Another technique suitable for 2-dimensional systems is to design the level set of a Lyapunov function as a \emph{polygon} whose boundary is crossed strictly inward, similarly to the comparison arguments \cite{tarabaStabilityLinearSecondorder2022}.

\emph{In this paper} we study a restricted version of the absolute stability problem: we consider a \emph{linear} feedback $u = \kappa(t) y$ and assume that not only its gain $\kappa$ is bounded, but the derivative $\dot\kappa$ is constrained as well. Our analysis is based on a slight modification of the variational approach. A comparison between the two is given in Remark~\ref{re: comparison to variational approach}.

\emph{Main results} (Section~\ref{se: results}) consist of the following:
\begin{itemize}
    \item A sufficient condition for absolute stability (Theorem~\ref{th: sufficient}) as well as a sufficient condition for the lack of absolute stability (Theorem~\ref{th: necessary}). These conditions are formulated in terms of the value function of an optimal control problem.
    \item A numerical algorithm for approximating the value function (Proposition~\ref{pr: numerical method}). In combination with the two theorems, it can be used as a test of stability.
\end{itemize}
Our stability condition is compared to the one found in~\cite{ignatyevStabilityLinearOscillator1997} by means of a practically motivated example in Section~\ref{se: example}.

Let us finish this introduction with a few comments:
\begin{itemize}
    \item For technical reasons discussed in Section~\ref{se: why linear feedback}, we consider linear feedback, unlike in the classical Lurie problem. This can be viewed as a linearization of the nonlinear Lurie system about the origin.
    \item The constraint on $\dot\kappa$ may depend on $\kappa$. This is different from many publications where the constraint essentially follows from the chosen analysis technique. Allowing a more general constraint is useful in some practical cases, such as the one treated in Section~\ref{se: example}.
    \item An easy extension of our results even to 3-dimensional systems is hindered by several factors. Firstly, we use polar coordinates and assume that the system is oscillatory. These concepts should be reexamined in a higher-dimensional situation. Secondly, our test of absolute stability is based on the numerical solution of a partial differential equation. As the number of dimensions grows, the numerical method becomes exponentially more expensive. The curse of dimensionality also concerns the \emph{multiple-input multiple-output} case with many independent feedback gains.
    \item In the classical 2-dimensional Lurie--Aizerman problem, the variational/comparison approach yields a \emph{necessary and sufficient} condition for absolute stability. The reason for why such a condition is missing from our work is explained in Section~\ref{se: gap}. However, in the extreme case of unconstrained $\dot\kappa$ our conditions converge to the necessary and sufficient conditions as the method itself relaxes into the variational/comparison method.
\end{itemize}

\section{Problem Statement}

Consider the two-dimensional SISO LTI system
\begin{subequations}
    \label{eq: system}
    \begin{alignat}{3}
        \dot x &= Ax + bu, \quad && x\in\bbR^2, \: u\in\bbR, \\
        y &= c^T x, && y\in\bbR
    \end{alignat}
    closed by the time-varying linear feedback
    \begin{equation}
        \label{eq: feedback}
        u = \kappa(t) y, \quad t \geq 0
    \end{equation}
    where the feedback gain $\kappa$ is a continuously differentiable function that satisfies
    \label{eq: feedback constraints}
    \begin{gather}
        \kappa(t) \in [0,1], \\
        \kdmin(\kappa(t)) \leq \dot\kappa(t) \leq \kdmax(\kappa(t))
    \end{gather}
\end{subequations}
for all $t \geq 0$. Here $\kdmin$ and $\kdmax$ are given continuous functions such that
\begin{equation}
    -\infty < \kdmin(\kappa) \leq 0 \leq \kdmax(\kappa) < \infty.
\end{equation}

The following is a straightforward adaptation of the classical absolute stability concept to our case (\emph{linear} feedback~\eqref{eq: feedback} and constrained gain derivative $\dot\kappa$).

\begin{definition}
    \label{def: absolute stability}
    System~\eqref{eq: system} is called \emph{absolutely stable} if its zero solution is globally asymptotically stable under every admissible feedback~\eqref{eq: feedback}.
\end{definition}

\begin{remark}
    The global character of asymptotic stability is actually irrelevant in our case since the closed loop is linear.
\end{remark}

\section{Preliminaries}

\subsection{Extended System}

Let us begin by attaching the feedback gain $\kappa$ to the system's state vector. This extends the system's dimensionality from 2 to 3. Furthermore, we expand the functional space for $\dot\kappa$ from continuous ($C^0$) to locally $L^1$-integrable ($\Loneloc$) functions. This move is motivated by the upcoming discussion of a variational problem that requires relaxed regularity. Since $C^0$ is $L^1$-dense in $\Loneloc$, the enlargement of the admissible space is negligible enough so that not to affect the absolute stability property. As a result, we obtain the extended closed-loop system
\begin{subequations}
    \label{eq: closed loop}
    \begin{align}
        \dot x &= (A + \kappa B) x, \quad B = bc^T, \\
        \dot \kappa &= \nu(\kappa, t), \\
        \kappa(0) &\in [0,1], \\
        \nu &\in \calA, \\
        \begin{split}
            \calA &= \setb{\nu \in \Loneloc :
            \kdmin(\kappa) \leq \nu \leq \kdmax(\kappa), \\
            &\phantom{ {}= \big\{} \nu \geq 0 \text{ if } \kappa = 0, \:
            \nu \leq 0 \text{ if } \kappa = 1}.
        \end{split}
    \end{align}
\end{subequations}

The uncertain parts of system~\eqref{eq: closed loop} are:
\begin{itemize}
    \item initial feedback gain $\kappa(0) \in [0,1]$;
    \item unspecified function $\nu\in\calA$.
\end{itemize}

\subsection{Change of Coordinates}

Here we transform the system to the polar coordinates and, assuming the oscillating character of the dynamics, use the polar angle as the new independent variable. The system's dimensionality drops back to 2. The original problem reduces to the partial stability of the new 2-dimensional system.

Let us introduce the basis vectors
\begin{equation}
    \bfi(\theta) = \begin{bmatrix}
        \cos\theta \\ \sin\theta
    \end{bmatrix}, \quad
    \bfj(\theta) = \begin{bmatrix}
        -\sin\theta \\ \phantom{-}\cos\theta
    \end{bmatrix}
\end{equation}
and define \emph{polar coordinates} $(r, \theta)$ via
\begin{equation}
    x = r \bfi(\theta), \quad r \geq 0.
\end{equation}

\begin{remark}
    \label{re: theta branch}
    Polar angle $\theta$ is a multivalued function of $x$ whose adjacent branches differ by $2\pi$. As time goes on, we require that $\theta$ changes \emph{continuously}, progressively moving from one branch to another as $x$ loops around the origin.
\end{remark}

In the polar coordinates, system \eqref{eq: closed loop} reads
\begin{subequations}
    \label{eq: system polar in time}
    \begin{align}
        \dot r &= h(\kappa, \theta) r, \\
        \dot\theta &= g(\kappa, \theta), \\
        \dot\kappa &= \nu(\kappa, t), \\
        \kappa(0) &\in [0,1], \\
        \nu &\in \calA
    \end{align}
\end{subequations}
where
\begin{subequations}
    \begin{align}
        h(\kappa, \theta) &= \bfi(\theta)^T (A + \kappa B) \bfi(\theta), \\
        g(\kappa, \theta) &= \bfj(\theta)^T (A + \kappa B) \bfi(\theta).
    \end{align}
\end{subequations}

Note that the qualitative behavior of system~\eqref{eq: system} may depend on $\kappa$: for different fixed values of $\kappa$ the system may be oscillating or non-oscillating (overdamped). In this paper we focus our attention on the case when the system remains oscillating for all admissible $\kappa$. Other cases require separate but similar treatment.

\begin{assumption}
    \label{as: theta increasing}
    System~\eqref{eq: system} is oscillatory for every choice of the feedback gain $\kappa \in [0,1]$ and oscillates in the counterclockwise direction. That is, in~\eqref{eq: system polar in time} we assume $\dot\theta = g(\kappa, \theta) > 0$ and thus, in view of Remark~\ref{re: theta branch}, $\theta \to \infty$ as $t \to \infty$.
\end{assumption}

\begin{remark}
    Assumption~\ref{as: theta increasing} can be ensured if matrix $A + \kappa B$ in \eqref{eq: closed loop} has a pair of nonreal eigenvalues for every $\kappa \in [0,1]$. Indeed, it implies that $g(\kappa, \theta)$ has constant sign. If it happens that $g < 0$ then its sign can be reversed by swapping the coordinates: take $x^\mathrm{new} = \begin{bmatrix}
        x_2 & x_1
    \end{bmatrix}^T$.
\end{remark}

With Assumption~\ref{as: theta increasing} established, let us pick $\theta$ as the new independent variable in place of time. It brings~\eqref{eq: system polar in time} to
\begin{subequations}
    \label{eq: system polar}
    \begin{align}
        r' &= f(\kappa, \theta) r, \\
        \kappa' &= g(\kappa, \theta) \nu(\kappa, \theta), \\
        \kappa(0) &\in [0,1], \\
        \nu &\in \calA
    \end{align}
\end{subequations}
where prime denotes the derivative with respect to $\theta$ and
\begin{equation}
    f(\kappa, \theta) = \frac{h(\kappa, \theta)}{g(\kappa, \theta)}.
\end{equation}
Note the abuse of notation in~\eqref{eq: system polar}: we write $\nu(\kappa, \theta)$ in place of $\nu(\kappa, t(\theta))$. It shall cause no confusion since the constraint $\nu \in \calA$ is time-independent.

Absolute stability of~\eqref{eq: system} translates to the absolute stability of~\eqref{eq: system polar} with respect to the radius $r$~-- the relevant term is ``partial stability'' \cite[Definition~2]{vorotnikovPartialStabilityControl2005}. Accordingly, we combine the concepts of absolute and partial stability in the following definition.

\begin{definition}
    \label{def: absolute r-stability}
    System~\eqref{eq: system polar} is called \emph{absolutely $r$-stable} if for every $\nu \in \calA$ the partial equilibrium $r = 0$ of~\eqref{eq: system polar} is globally asymptotically $r$-stable, i.e.:
    \begin{itemize}
        \item for every $\theta_0 \geq 0$ and arbitrary initial conditions $r(\theta_0) \geq 0$, $\kappa(\theta_0) \in [0,1]$ the solution $(r(\theta), \kappa(\theta))$ exists for all $\theta \geq \theta_0$;
        \item for each $\theta_0 \geq 0$ and $\epsilon > 0$ there exists $\delta > 0$ such that $r(\theta_0) < \delta$ implies $r(\theta) < \epsilon$ for all $\theta \geq \theta_0$ and arbitrary $\kappa(\theta_0) \in [0,1]$;
        \item $r(\theta) \to 0$ as $\theta \to \infty$ for every solution $(r(\theta), \kappa(\theta))$.
    \end{itemize}
\end{definition}

The following obvious statement is made for the sake of clarity.

\begin{proposition}
    \label{pr: equivalence}
    System~\eqref{eq: system} is absolutely stable if and only if~\eqref{eq: system polar} is absolutely $r$-stable.
\end{proposition}

\subsection{Variational Approach}

Let us formulate an \emph{optimal control problem} for system~\eqref{eq: system polar} on a fixed interval $[0, \bar\theta]$: given $\bar\kappa \in [0,1]$, find a ``control'' $\nu\in\calA$ such that
\begin{equation}
    \label{eq: optimization problem}
    \frac{r(\bar\theta)}{r(0)} \to \max \text{ subject to } \kappa(\bar\theta) = \bar\kappa.
\end{equation}
The initial values $r(0)$ and $\kappa(0)$ are free.

\begin{remark}
    \label{re: comparison to variational approach}
    Our main results in Section~\ref{se: results} are delivered by a solution to the optimization problem~\eqref{eq: optimization problem}. This is somewhat similar to the classical variational analysis of~\cite{pyatnitskiiAbsoluteStabilityNonstationary1970} and \cite{margaliotStabilityAnalysisSwitched2006} with some differences. In~\cite{margaliotStabilityAnalysisSwitched2006}, the full norm of the endpoint is maximized. We, on the other hand, maximize only the $r$-coordinate while the $\kappa$-coordinate is treated as a parameter since we are investigating partial $r$-stability. In~\cite{margaliotStabilityAnalysisSwitched2006}, the initial point is fixed; later, it is shown that the value function of the optimization problem is homogeneous with respect to the initial point. We leave the initial point free and utilize $r$-homogeneity directly in the optimization problem~\eqref{eq: optimization problem} via normalizing the cost by $r(0)$. This helps to keep the number of variables down and facilitates a numerical solution.
\end{remark}

We now aim to derive the \emph{Hamilton--Jacobi--Bellman} (HJB) equation for problem~\eqref{eq: optimization problem}, i.e., a partial differential equation for the \emph{value function} $\rho(\kappa, \theta)$ defined as
\begin{equation}
    \label{eq: value definition}
    \rho(\bar\kappa, \bar\theta) = \max\frac{r(\bar\theta)}{r(0)} \text{ subject to } \kappa(\bar\theta) = \bar\kappa.
\end{equation}

\begin{remark}
    The value function $\rho$ is related to the absolute Lyapunov exponent from~\cite{barabanovAbsoluteCharacteristicExponent1989}: indeed,
    \begin{equation}
        \limsup_{\theta\to\infty} \max_{\kappa\in[0,1]} \frac{\ln\rho(\kappa,\theta)}{\theta}
    \end{equation}
    may be called the ``partial'' absolute Lyapunov exponent~-- that is, exponent only with respect to $r$.
\end{remark}

Proceeding with the derivation of the HJB equation, consider two angles $\theta_{1,2}$ such that $0 \leq \theta_1 < \theta_2$. There exists a trajectory $(r, \kappa)$ of~\eqref{eq: system polar} on the interval $[0, \theta_1]$ with an arbitrarily prescribed endpoint that attains the maximum in~\eqref{eq: value definition}:
\begin{equation}
    \label{eq: rho at theta1}
    \rho\pb{\kappa(\theta_1), \theta_1} = \frac{r(\theta_1)}{r(0)}.
\end{equation}
The value of $r(\theta_1)$ can be viewed as arbitrary because system~\eqref{eq: system polar} as well as the ratio $r(\theta_1)/r(0)$ are invariant with respect to the scaling of $r$. Let us extend this trajectory to the interval $[\theta_1, \theta_2]$ arbitrarily. Then due to~\eqref{eq: value definition}
\begin{equation}
    \label{eq: rho at theta2}
    \rho\pb{\kappa(\theta_2), \theta_2} \geq \frac{r(\theta_2)}{r(0)}.
\end{equation}
From~\eqref{eq: rho at theta1} and~\eqref{eq: rho at theta2}
\begin{equation}
    \frac{r(\theta_2)}{\rho\pb{\kappa(\theta_2), \theta_2}}
    \leq
    \frac{r(\theta_1)}{\rho\pb{\kappa(\theta_1), \theta_1}}.
\end{equation}
Since $(r, \kappa)$ is a free trajectory on $[\theta_1, \theta_2]$, we conclude that the ratio
\begin{equation}
    \label{eq: ratio}
    \gamma(r, \kappa, \theta) = \frac{r}{\rho(\kappa, \theta)}
\end{equation}
is non-increasing along all solutions of~\eqref{eq: system polar}.

Let us assume for a moment that $\rho$ is smooth. Then the derivative of $\gamma$ along the solutions of~\eqref{eq: system polar} must be nonpositive:
\begin{equation}
    \label{eq: negative derivative}
    \frac{\dif\gamma}{\dif\theta} = -\frac{r}{\rho^2} \pb{
        \rho_\theta
        + g(\kappa, \theta) \nu(\kappa, \theta) \rho_\kappa
        - f(\kappa, \theta) \rho
    } \leq 0
\end{equation}
(subscripts here represent partial derivatives). Condition~\eqref{eq: negative derivative} is equivalent to
\begin{equation}
    \label{eq: derivative inequality}
    \rho_\theta
    + g(\kappa, \theta) \nu(\kappa, \theta) \rho_\kappa
    \geq f(\kappa, \theta) \rho.
\end{equation}
For~\eqref{eq: derivative inequality}, it is sufficient to require
\begin{equation}
    \label{eq: HJB equation}
    \rho_\theta
    + g(\kappa, \theta) \min_{\nu \in \calA} \setb{\nu(\kappa, \theta) \rho_\kappa}
    = f(\kappa, \theta) \rho
\end{equation}
because $g > 0$ by Assumption~\ref{as: theta increasing}. Equation~\eqref{eq: HJB equation} is a partial differential equation describing a function $\rho$ of two independent variables $\theta$ and $\kappa$.

\begin{figure}[t]
    \centering
    \includegraphics{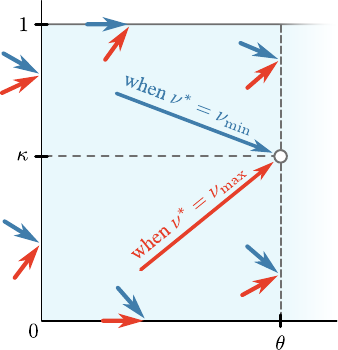}
    \caption{Possible characteristic directions of the HJB equation~\eqref{eq: HJB equation} depending on the value $\nu^*$ of $\nu$ that attains a minimum in~\eqref{eq: HJB equation}.}
    \label{fig: characteristics}
\end{figure}

We claim that~\eqref{eq: HJB equation} is the desired HJB equation. Before we prove it formally, let us discuss the \emph{boundary conditions} the equation needs. The domain of~\eqref{eq: HJB equation} is the half-strip $[0, \infty) \times [0,1]$ in the $(\theta, \kappa)$ plane. To determine the necessary boundary conditions, we look at the characteristic directions. Depending on the sign of the partial derivative $\rho_\kappa(\kappa, \theta)$, the characteristic direction \emph{entering} a point $(\kappa, \theta)$ is
\begin{equation}
    \label{eq: characteristics}
    \text{either }
    \begin{bmatrix}
        1 \\
        g(\kappa, \theta) \numax(\kappa)
    \end{bmatrix}
    \text{ or }
    \begin{bmatrix}
        1 \\
        g(\kappa, \theta) \numin(\kappa)
    \end{bmatrix}
\end{equation}
where
\begin{subequations}
    \begin{align}
        \numin(\kappa) &= \begin{cases}
            \kdmin(\kappa), \,& \kappa \in [0,1), \\
            0, & \kappa = 1,
        \end{cases} \\
        \numax(\kappa) &= \begin{cases}
            0, & \kappa = 0, \\
            \kdmax(\kappa), & \kappa \in (0,1].
        \end{cases}
    \end{align}
\end{subequations}
These two options are shown in Fig.~\ref{fig: characteristics}. A characteristic direction is the direction of the solution propagation and can be mechanically interpreted as the direction of the flow. A boundary condition is required at the points where the flow is coming in from the outside of the domain. Such are the boundary points with $\theta = 0$. Therefore, equation~\eqref{eq: HJB equation} requires a boundary condition for $\theta = 0$. Plugging $\theta = 0$ into~\eqref{eq: value definition} we obtain
\begin{equation}
    \label{eq: boundary condition}
    \rho(\kappa, 0) = 1 \text{ for all } \kappa\in[0,1].
\end{equation}
On the rest of the boundary ($\kappa = 0$ and $\kappa = 1$) the flow is not entering the domain; no boundary condition is needed there.

As for the notion of a solution of the HJB equation, we adopt the established concept of a \emph{viscosity solution}~\cite[Sec.~5.3]{liberzonCalculusVariationsOptimal2012}. The following lemma confirms that~\eqref{eq: HJB equation} is indeed the HJB equation in the sense that the value function~\eqref{eq: value definition} is a viscosity solution of~\eqref{eq: HJB equation}.

\begin{lemma}
    \label{le: HJB}
    Equation~\eqref{eq: HJB equation} is the Hamilton--Jacobi--Bellman equation for the variational problem~\eqref{eq: optimization problem}. That is, value function~\eqref{eq: value definition} is the viscosity solution of~\eqref{eq: HJB equation} with the boundary condition~\eqref{eq: boundary condition}.
\end{lemma}

\begin{proof}
    The proof essentially follows the discussion in~\cite[Sec.~5.3.3]{liberzonCalculusVariationsOptimal2012}. A viscosity solution $\rho$ is, in particular, a viscosity \emph{super}solution which means the following. Consider a smooth function $\rho\sub(\kappa, \theta)$, called a \emph{test function} in~\cite{liberzonCalculusVariationsOptimal2012}, that satisfies $\rho\sub = \rho$ at a point $(\kappa_0, \theta_0)$ and $\rho\sub \leq \rho$ in a neighborhood of that point. Then equation~\eqref{eq: HJB equation} must hold for $\rho\sub$ in the neighborhood of $(\kappa_0, \theta_0)$ as a ``$\geq$'' inequality:
    \begin{equation}
        \rho\sub_\theta
        + g(\kappa, \theta) \min_{\nu \in \calA} \setb{\nu(\kappa, \theta) \rho\sub_\kappa}
        \geq f(\kappa, \theta) \rho\sub.
    \end{equation}
    If we formed the ratio $\gamma$ in~\eqref{eq: ratio} with $\rho\sub$ instead of $\rho$ then such $\gamma$ would be smooth and, from the latter inequality, non-increasing along the solutions of~\eqref{eq: system polar}:
    \begin{equation}
        \frac{\dif}{\dif\theta} \px{\frac{r}{\rho\sub(\kappa, \theta)}} \leq 0
    \end{equation}
    at the point $(\kappa_0, \theta_0)$ for all $r \geq 0$. Therefore, since in the neighborhood of $(\kappa_0, \theta_0)$
    \begin{equation}
        \gamma(r, \kappa, \theta) = \frac{r}{\rho(\kappa, \theta)}
        \leq \frac{r}{\rho\sub(\kappa, \theta)}
    \end{equation}
    we conclude that $\gamma$ defined as~\eqref{eq: ratio} is also non-increasing along the solutions of~\eqref{eq: system polar}.  Traversing the chain of implications backwards from~\eqref{eq: ratio} to~\eqref{eq: value definition}, we arrive at
    \begin{equation}
        \rho(\kappa, \theta) \geq \max\frac{r(\theta)}{r(0)}
    \end{equation}
    along every solution of~\eqref{eq: system polar}. To finish the proof, one would similarly use the fact that $\rho$ is a viscosity \emph{sub}solution to show that the latter inequality is in fact an equality.
\end{proof}

\section{Main Results}
\label{se: results}

\subsection{Conditions for Stability and Instability}

The following results are based on the variational property~\eqref{eq: value definition}.

\begin{theorem}[Stability]
    \label{th: sufficient}
    Under Assumption~\ref{as: theta increasing}, let $\rho$ be the viscosity solution of the HJB equation~\eqref{eq: HJB equation} with boundary condition~\eqref{eq: boundary condition}. If
    \begin{equation}
        \rho(N\pi, \kappa) < 1 \text{ for some } N \in \bbN \text{ and all } \kappa\in[0,1]
    \end{equation}
    then system~\eqref{eq: system} is absolutely stable per Definition~\ref{def: absolute stability}.
\end{theorem}

\begin{proof}
    By Lemma~\ref{le: HJB}, $\rho$ satisfies~\eqref{eq: value definition} along every solution of~\eqref{eq: system polar}. Therefore, $r(N\pi) < r(0)$ for every solution, and so the map $r(0) \mapsto r(N\pi)$ is contractive. Since system~\eqref{eq: system polar} is $\pi$-periodic with respect to $\theta$, this proves that $r$ converges to zero exponentially. Thus, \eqref{eq: system polar} is absolutely $r$-stable. By Proposition~\ref{pr: equivalence} system~\eqref{eq: system} is absolutely stable as well.
\end{proof}

\begin{theorem}[Instability]
    \label{th: necessary}
    If in Theorem~\ref{th: sufficient}
    \begin{equation}
        \rho(N\pi, \kappa) > 1 \text{ for some } N \in \bbN \text{ and all } \kappa\in[0,1]
    \end{equation}
    then system~\eqref{eq: system} is \emph{not} absolutely stable per Definition~\ref{def: absolute stability}.
\end{theorem}

\begin{proof}
    Consider the maximizer $\nu\in\calA$ of the variational problem~\eqref{eq: optimization problem} on the interval $[0, N\pi]$. This function $\nu$ attains the maximum of $r(N\pi)/r(0)$~-- that is, $\nu$ is the ``most destabilizing'' feedback in the terminology of~\cite{margaliotStabilityAnalysisSwitched2006}. This $\nu$ destabilizes the system because with it $r(N\pi) > r(0)$ independently of $\kappa(0)$. Furthermore, $\nu$ can be approximated arbitrarily closely in the $L^1$-norm by a continuous function $\dot\kappa$ whose integral $\kappa$ is an admissible feedback gain. This $\kappa$ destabilizes the original system~\eqref{eq: system}. Therefore, \eqref{eq: system} is not absolutely stable.
\end{proof}

\subsection{Numerical Method}

We supplement Theorems~\ref{th: sufficient} and~\ref{th: necessary} with a numerical scheme to solve the HJB equation~\eqref{eq: HJB equation}. The scheme is of the semi-Lagrangian type \cite{falconeSemiLagrangianApproximationSchemes2014}. It combines the Eulerian point of view (approximation on a fixed grid) and the Lagrangian one (propagation of the solution along the characteristics).

Consider a uniform grid in the $(\kappa, \theta)$ domain with steps $\dlt\kappa$ and $\dlt\theta$ satisfying the following assumption.

\begin{assumption}
    \label{as: CFL}
    The Courant--Friedrichs--Lewy (CFL) condition holds:
    \begin{equation}
        \dlt\theta < \frac{\dlt\kappa}{
            \max\limits_{\substack{
                \nu\in\calA \\
                \kappa \in [0,1] \\
                \theta \geq 0
            }} \absb{g(\kappa, \theta) \nu(\kappa, \theta)}
        }.
    \end{equation}
\end{assumption}

\begin{remark}
    \label{re: grid size}
    Assumption~\ref{as: CFL} is ensured by choosing
    \begin{equation}
        \dlt\theta = \frac{\dlt\kappa}{
            (\norm{A} + \norm{B})
            \max\limits_{\kappa \in [0,1]} \setb{\abs{\numin(\kappa)}, \abs{\numax(\kappa)}}
        }.
    \end{equation}
\end{remark}

\begin{remark}
    Although the CFL condition is not strictly necessary for the numerical stability of semi-Lagrangian methods, we still assume that it holds as it makes it easier to ensure that characteristics do not violate the bound $\kappa \in [0,1]$.
\end{remark}

Approximation of $\rho$ on the grid is given by the values
\begin{equation}
    \rho_j^n \approx \rho(\kappa_j, \theta_n), \quad
    \kappa_j = j\dlt\kappa, \quad
    \theta_n = n\dlt\theta.
\end{equation}
We also define $\tilde\rho(\kappa, \theta_n)$ for an arbitrary $\kappa \in [0,1]$ as a linear interpolant between the nodes of the grid:
\begin{multline}
    \tilde\rho(\kappa, \theta_n) = \rho_j^n
    + \frac{\kappa - \kappa_j}{\dlt\kappa}
    (\rho_{j+1}^n - \rho_j^n) \\
    \text{ for } \kappa \in [\kappa_j, \kappa_{j+1}).
\end{multline}

To implement the Lagrangian part of the algorithm, we look again at two possible characteristic directions~\eqref{eq: characteristics}. The sign of the partial derivative $\rho_\kappa$ determines which characteristic is the actual one at a given point. Since the sign is not known beforehand, we consider both of these directions. Discarding the nonlinear terms, admissible characteristics~\eqref{eq: characteristics} yield two options for the propagation rules on the Eulerian grid: respectively,
\begin{subequations}
    \label{eq: propagation options}
    \begin{equation}
        \rho_j^{n+1}
        = \tilde\rho\pb{\kappa_j
        - g_j^n \numax(\kappa_j) \dlt\theta,
        \theta_n}
        + f_j^n \rho_j^n \dlt\theta
    \end{equation}
    and
    \begin{equation}
        \rho_j^{n+1}
        = \tilde\rho\pb{\kappa_j
        - g_j^n \numin(\kappa_j) \dlt\theta,
        \theta_n}
        + f_j^n \rho_j^n \dlt\theta
    \end{equation}
\end{subequations}
where $f_j^n = f(\kappa_j, \theta_n)$ and $g_j^n = g(\kappa_j, \theta_n)$. Formulas~\eqref{eq: propagation options} use the interpolant $\tilde\rho$ to approximate the value of $\rho$ at the point $(\dots, \theta_n)$ preceding $(\kappa_j, \theta_{n+1})$ along a characteristic: generally, the predecessor falls between the grid nodes. Due to Assumption~\ref{as: CFL}, the $\kappa$-coordinate of the predecessor is no further than one $\dlt\kappa$ step away from $\kappa_j$; therefore, it remains within $[0,1]$.

It can be seen directly from~\eqref{eq: HJB equation} or from the variational characterization~\eqref{eq: optimization problem} of $\rho$ that the choice between the two options~\eqref{eq: propagation options} should be in favor of one that produces the larger $\rho_j^{n+1}$. Thus, we attain our algorithm.

\begin{proposition}
    \label{pr: numerical method}
    Under Assumptions~\ref{as: theta increasing} and \ref{as: CFL}, the viscosity solution of the HJB equation~\eqref{eq: HJB equation} with boundary condition~\eqref{eq: boundary condition} is approximated by the finite difference scheme
    \begin{align}
        \rho_j^{n+1} = \max \setB{
            &\tilde\rho\pb{\kappa_j
            - g_j^n \numax(\kappa_j) \dlt\theta,
            \theta_n}, \notag \\
            &\tilde\rho\pb{\kappa_j
            - g_j^n \numin(\kappa_j) \dlt\theta,
            \theta_n}
        } + f_j^n \rho_j^n \dlt\theta
    \end{align}
    for $n = 0,1,\dots$ with initial condition $\rho_j^0 = 1$.
\end{proposition}

\section{Discussion}

Let us discuss a couple of limitations of this work.

\subsection{Why Linear Feedback?}
\label{se: why linear feedback}

In contrast to the classical Lurie problem with nonlinear feedback $u = \phi(y,t)$, we only consider $u = \kappa(t) y$. In many publications on absolute stability the nonlinear case is implicitly reduced to the linear one since along every particular trajectory $\kappa(t)$ can be found such that $\phi(y,t) = \kappa(t) y$. However, once we want to introduce derivative constraints, it is necessary to consider
\begin{equation}
    \dot\kappa(t) = \frac{\dif}{\dif t} \px{\frac{\phi(y,t)}{y}}
    = \frac{\phi_t y + \phi_y y \dot y - \phi \dot y}{y^2}.
\end{equation}
From here, one can derive constraints on $\dot\kappa$ assuming that the partial derivatives $\phi_t$ and $\phi_y$ are somehow constrained. The following points complicate the analysis, however:
\begin{enumerate}
    \item Derivative $\dot\kappa$ is unbounded whenever the output $y$ crosses zero. Our numerical method (Proposition~\ref{pr: numerical method}) has to be modified since currently it assumes the CFL condition (Assumption~\ref{as: CFL}) which implies that the admissible range of $\nu = \dot\kappa$ is bounded.
    \item Reduction to the linear feedback is conservative in presence of the derivative constraints. For example, consider the special case of a constraint $\phi_t = 0$, i.e., $u = \phi(y)$. This implies a nonlocal constraint on the linear gain $\kappa$: $\kappa(t_1) = \kappa(t_2)$ whenever $y(t_1) = y(t_2)$. Such a constraint cannot be accurately captured by a condition on $\dot\kappa$. Thus, it is not clear how strong of a result this approach can yield.
    \item The constraint on $\dot\kappa$ may depend on $\dot y$ and therefore on $x$ which makes the closed-loop system nonlinear with respect to $x$. Unlike in the Lurie--Aizerman case, it is possible to have local but not global asymptotic stability. If we do not wish to resort to local linearization-based analysis, then a discussion of the domain of attraction is required.
\end{enumerate}

\subsection{The Necessity~-- Sufficiency Gap}
\label{se: gap}

There is a visible gap between the conditions of Theorems~\ref{th: sufficient} and \ref{th: necessary}: if $\rho(N\pi, \kappa)$ is above or below 1 depending on $\kappa$ then we cannot make a decision about absolute stability. One can only check if a Theorem works for a larger $N$.

It is difficult to make any conclusion in the intermediate case because it allows both $r(N\pi) < r(0)$ and $r(N\pi) > r(0)$ depending on $\kappa(0)$. In our variational problem~\eqref{eq: optimization problem}, we discard the information about $\kappa(0)$ treating it as one of the optimizable degrees of freedom. It might be possible to resolve the issue by tracing the optimal trajectories from the starting to the endpoint. However, attention should be paid to the possibility of sliding dynamics induced by the switching behavior of the optimal control.

\section{Example}
\label{se: example}

In~\cite{ponomarevNonlinearAnalysisSynchronous2024} we came across an equation of the form
\begin{equation}
    \label{eq: example equation}
    m(t) \ddot\xi + k \dot\xi + \xi = 0
\end{equation}
where $k > 0$ and the uncertain function $m$ satisfies
\begin{equation}
    \label{eq: example constraint}
    0 < \mmin \leq m(t) \leq \mmax < \infty.
\end{equation}
Equation~\eqref{eq: example equation} appears during the linearized \emph{incremental stability} analysis of a synchronization circuit used in power grids. In an almost perfectly balanced grid, $m$ is almost constant: $\mmin \approx \mmax$. Let us now assume, in addition to~\eqref{eq: example constraint}, that
\begin{equation}
    \label{eq: example derivative constraint}
    \absb{\dot m(t)} \leq \epsilon.
\end{equation}
This also makes sense in a nearly balanced grid.

We represent~\eqref{eq: example equation} as a Lurie system by introducing the coordinates $x_1 = \dot\xi$ and $x_2 = \xi$. Then
\begin{subequations}
    \label{eq: example system final form}
    \begin{align}
        \dot x_1 &= -\frac{1}{m(t)} (kx_1 + x_2), \\
        \dot x_2 &= x_1.
    \end{align}
\end{subequations}
Let us represent
\begin{equation}
    -\frac{1}{m(t)} = -\frac{1}{\mmin} + \kappa(t) \px{\frac{1}{\mmin} - \frac{1}{\mmax}}
\end{equation}
so that $\kappa(t) \in [0,1]$. Then $x$ is governed by~\eqref{eq: system} with
\begin{multline}
    A = \begin{bmatrix}
        -k/\mmin & -1/\mmin \\ 1 & 0
    \end{bmatrix}, \\
    b = \px{\frac{1}{\mmin} - \frac{1}{\mmax}} \begin{bmatrix}
        1 \\ 0
    \end{bmatrix}, \quad c = \begin{bmatrix}
        k \\ 1
    \end{bmatrix}.
\end{multline}
To determine the constraints on $\dot\kappa$, we derive
\begin{equation}
    \dot\kappa(t) = -\dot m(t) \px{\kappa(t) + \frac{\mmin}{\mmax - \mmin}}
\end{equation}
and due to~\eqref{eq: example derivative constraint} obtain
\begin{equation}
    \abs{\dot\kappa} \leq \epsilon \px{\kappa + \frac{\mmin}{\mmax - \mmin}}.
\end{equation}
Thus,
\begin{subequations}
    \begin{align}
        \kdmin(\kappa) &= -\epsilon \px{\kappa + \frac{\mmin}{\mmax - \mmin}}, \\
        \kdmax(\kappa) &= \phantom{ - } \epsilon \px{\kappa + \frac{\mmin}{\mmax - \mmin}}.
    \end{align}
\end{subequations}

\begin{figure}[t]
    \centering
    \includegraphics{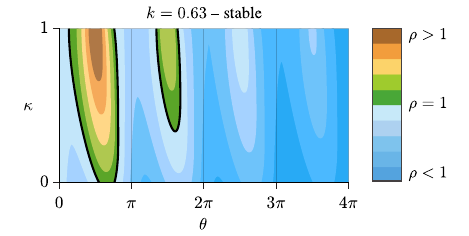}

    \vspace{2mm}
    \includegraphics{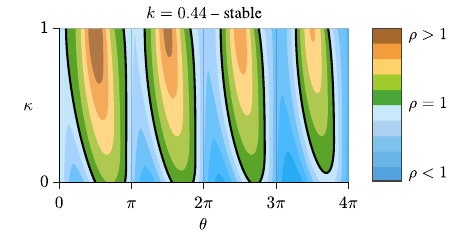}
    
    \vspace{2mm}
    \includegraphics{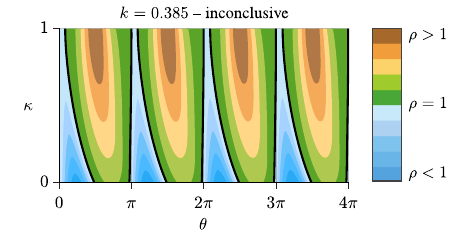}
    
    \vspace{2mm}
    \includegraphics{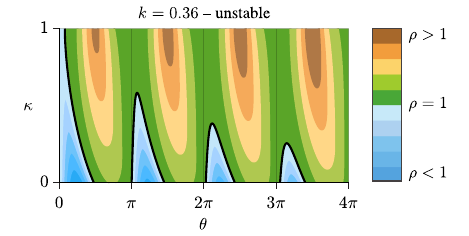}

    \caption{Approximate solution $\rho(\kappa, \theta)$ of the HJB equation~\eqref{eq: HJB equation} in the example, computed for different values of $k$ by the method of Proposition~\ref{pr: numerical method}. By Theorems~\ref{th: sufficient} and \ref{th: necessary}, absolute stability is lost as $k$ goes below about 0.385.}
    \label{fig: example}
\end{figure}

Suppose, for this example, that $\mmin = 1$ and $\mmax = 4$. According to~\cite[Lemma~2]{ponomarevNonlinearAnalysisSynchronous2024}, equation~\eqref{eq: example equation} is exponentially stable if
\begin{equation}
    \label{eq: example old condition}
    k > \sqrt{\mmax} - \sqrt{\mmin} = 1.
\end{equation}
For such $k$, we had absolute stability without restricting $\dot m$. Note that this condition is conservative. It is attained by searching for a time-invariant quadratic Lyapunov function and, therefore, is the same as would have been found from the circle criterion. The comparison technique~\cite{leonovNecessarySufficientConditions2005}, on the other hand, yields a bound
\begin{equation}
    k \gtrsim 0.63
\end{equation}
which is, by construction, close to the exact threshold.

Let us now pick $\epsilon = 1$ and see how far down below 0.63 we can push $k$ while preserving absolute stability under the constraint~\eqref{eq: example derivative constraint}. We apply the algorithm of Proposition~\ref{pr: numerical method} to approximate the solution $\rho$ of the HJB equation~\eqref{eq: HJB equation} with grid steps $\dlt\kappa = 0.01$ and $\dlt\theta$ provided in Remark~\ref{re: grid size}. Then we conclude absolute stability or instability from Theorems~\ref{th: sufficient} and \ref{th: necessary}. The results are shown in Fig.~\ref{fig: example}: absolute stability is preserved while
\begin{equation}
    \label{eq: example result}
    k \gtrsim 0.385.
\end{equation}

Finally, let us compare this result to~\cite{ignatyevStabilityLinearOscillator1997} where a condition is given for the asymptotic stability of a second-order linear time-varying equation. Equation~\eqref{eq: example equation} can be written as
\begin{equation}
    \label{eq: example divided}
    \ddot\xi + \frac{k}{m(t)} \dot\xi + \frac{1}{m(t)} \xi = 0.
\end{equation}
By~\cite[Theorem~1]{ignatyevStabilityLinearOscillator1997}, \eqref{eq: example divided} is asymptotically stable if
\begin{equation}
    \frac{1}{2} \frac{\frac{\dif}{\dif t}(1/m)}{1/m} + \frac{k}{m} > \mathrm{const} > 0.
\end{equation}
A sufficient condition is thus
\begin{equation}
    k > \frac{\dot m}{2} \impliedby k > \frac{\epsilon}{2} = 0.5
\end{equation}
which is slightly more conservative than~\eqref{eq: example result}.

\section{Conclusion}

We have applied the variational method of absolute stability analysis to a system with constrained rate of change of the feedback gain. Using the Hamilton--Jacobi--Bellman theory, we analyze stability based on the value function of an optimal control problem (Theorems~\ref{th: sufficient} and \ref{th: necessary}). A practical implementation of the method relies on a numerical approximation of the value function (Proposition~\ref{pr: numerical method}). The numerical algorithm is easily implementable. With it, we have strengthened the previously known conditions for absolute stability in a practically motivated example.

The proof of the stability result (Theorem~\ref{th: sufficient}) can be straightforwardly extended to derive an exponential upper bound on the solutions of the closed-loop system. Furthermore, similarly to the techniques developed in~\cite{musaevaConstructionInvariantLyapunov2023,tarabaStabilityLinearSecondorder2022}, one may be able to characterize absolute stability further by constructing a Lyapunov function.

\bibliographystyle{IEEEtran}
\bibliography{lit}

\end{document}